\documentclass[12pt]{article}
\usepackage{amsmath}
\usepackage{amsthm}
\usepackage{amssymb}
\usepackage{amsfonts}
\usepackage{amscd}
\theoremstyle{plain}
	\newtheorem{thm}{Theorem}[section]
	\newtheorem*{cor}{Corollary}
	\newtheorem{lem}[thm]{Lemma}
	\newtheorem{prop}[thm]{Proposition}
\theoremstyle{definition}
	\newtheorem{defn}{Definition}[section]

\theoremstyle{remark}
	\newtheorem*{rem}{Remark}
	\newtheorem*{note}{Note}
	
	\newtheorem*{pf}{Proof}
\def\C{{\mathbb C}}
\def\Z{{\mathbb Z}}

\begin{document}
\title{K.~Saito's Duality for Regular Weight Systems \\
	and \\
	Duality for Orbifoldized Poincare Polynomials}
\author{Atsushi Takahashi
	\thanks{{\it E-mail address}~: atsushi@kurims.kyoto-u.ac.jp} \\
	\\
	Research Institute for Mathematical Sciences,\\
	Kyoto University,\\
	Kyoto 606-01, Japan}
\date{}
\maketitle
\begin{abstract}
We will show that the duality for regular weight system introduced by
K.~Saito can be interpreted as the duality for the orbifoldized
Poincare polynomial $\chi(W,G)(y,\bar{y})$.
\end{abstract}
%
\section*{Introduction}
~~~In \cite{ar:1}, Arnold discovered a strange duality among the 14
exceptional singularities.
This was interpreted by Dolgachev, Nikulin and Pinkham in terms of the duality
between algebraic cycles and transcendental cycles on certain K3 surfaces
\cite{dn:1}\cite{pi:1}.
Recently, K.~Saito discovered a new duality for regular weight
systems which contains the self-duality of ADE and Arnold's strange duality.
His theory of regular weight systems was originally developed in order to
understand the flat structure in the period map for primitive forms
\cite{sa:1}.
The theory of primitive forms can be interpreted as topological
Landau--Ginzburg models coupled to gravity \cite{ta:1}.
On the other hand, some LG (orbifold) models \cite{iv:1} have a duality
property in the sense \cite{bh:1}\cite{bh:2}\cite{ka:1}.
This fact is the motivation of this paper.
The organization of this paper is as follows.
In section 1, we prepare some definitions of the weight systems and
introduce the notion of P-duality for regular weight systems which
is a duality for the orbifoldized Poincare polynomials.
In section 2, we review the K.~Saito's theory of regular weight systems
and introduce the notion of M-duality for regular weight systems.
In section 3, we prove that M-duality is equivalent to P-duality.
\subsection*{Acknowledgement}
~~~I am deeply grateful to Professor K.~Saito for his encouragement.
I also would like to thank Professor S.~Hosono and Professor Toshiya Kawai
for useful discussions.
%
\newpage
\section{P-Duality for Regular Weight Systems}
~~~
\begin{defn}
We call $W:=(a_1,\dots ,a_n;h)$ a {\it weight system} if $a_i$, $i=1,\dots ,n$
and $h$ are positive integers such that $\max(a_i)< h$.
\end{defn}
The integers $a_i$ are called {\it weights} of $W$ and $h$ is called
the {\it coxeter number} of  $W$.
We also assume that $a_i\le h/2$ for all $i$.
\begin{defn}
$W$ is called {\it reduced} if $\gcd(a_1,\dots ,a_n;h)=1$.
\end{defn}
\begin{defn}~\cite{sa:2}
A weight system $W$ is called {\it regular} if a rational function:
\begin{equation}
\chi_W(T):=T^{\epsilon_W}\prod_{i=1}^{n}\frac{1-T^{h-a_i}}{1-T^{a_i}},
\end{equation}
has poles at most at $T=0$ where $\epsilon_W:=\sum_{i=1}^na_i-h$.
\end{defn}
We will treat only reduced regular weight systems which correspond to the
isolated hypersurface singularity with a $\C^*$-action.
We take a weighted homogeneous polynomial $F$ which has an isolated critical
point at $0$ and
\begin{equation}
F(\lambda^{a_1}z_1,\dots ,\lambda^{a_n}z_n)=\lambda^h F(z_1,\dots ,z_n),
~~~(a_i,h)=1,~~i=1,\dots ,n,
\end{equation}
then define $W=(a_1,\dots ,a_n;h)$.
For such a reduced regular weight system, there exist finite number of
integers $m_1<m_2\le\dots\le m_{\mu_W-1}<m_{\mu_W}$ called {\it exponents}
such that
\begin{equation}
\chi_W(T)=T^{m_1}+\dots+T^{m_{\mu_W}},
\end{equation}
where
\begin{equation}
\mu_W:=\prod_{i=1}^n\frac{h-a_i}{a_i}
\end{equation}
is called the {\it rank} of the weight system.
\begin{prop}
\begin{equation}
m_i+m_{\mu_W-i+1}=h,~~~for~~~i=1,\dots ,\mu_W.
\end{equation}
\end{prop}
Let $G$ be a discrete subgroup of $GL(n,\C)$ acting on
$(z_1,\dots ,z_n)$ such that $F$ is invariant under $G$.
We will also assume that $G$ is abelian and their elements are of the form
${\rm diag}({\bf e}[\omega_1\alpha_1],\dots ,{\bf e}[\omega_n\alpha_n])$ where
${\bf e}[\cdot]:=\exp(2\pi\sqrt{-1}\cdot)$, $\omega_i:=a_i/h$ and
$\alpha_i\in\Z$.
We call the group generated by
${\rm diag}({\bf e}[\omega_1],\dots ,{\bf e}[\omega_n])$
the {\it principal discrete group} and denote it by $G_0$.
We define the {\it Poincare polynomial}
$\chi(W,G)(y,\bar{y})$ for a pair$(W,G)$ as follows:
\begin{defn}[~\cite{iv:1}~\cite{ka:1} Vafa's fomula]
\begin{equation}
\chi(W,G)(y,\bar{y}):=\frac{(-1)^{n}}{|G|}\sum_{\alpha\in G}\chi_\alpha(W,G)
(y,\bar{y}),
\end{equation}
\begin{equation}
\begin{split}
\chi_\alpha(W,G)(y,\bar{y}):=& \sum_{\beta\in G}
\prod_{\omega_i\alpha_i\not\in\Z}\left(y\bar{y}\right)^{\frac{1-2\omega_i}{2}}
\left(\frac{y}{\bar{y}}\right)
^{-\omega_i\alpha_i+\left[\omega_i\alpha_i\right]+\frac{1}{2}}\\
& \times\prod_{\omega_i\alpha_i\in\Z}{\bf e}\left[\omega_i\beta_i+\frac{1}{2}
\right]\frac{1-{\bf e}\left[(1-\omega_i\beta_i)\right](y\bar{y})^{1-\omega_i}}
{1-{\bf e}\left[\omega_i\beta_i\right](y\bar{y})^{\omega_i}},
\end{split}
\end{equation}
where $[\omega_i\alpha_i]$ denotes the greatest integer smaller than
$\omega_i\alpha_i$.
\end{defn}
\begin{rem}
Note that if we put $T^h=y\bar{y}$, we have
\begin{equation}
\chi_W(T)=(y\bar{y})^{\frac{\epsilon_W}{h}}\chi(W)(y,\bar{y}),
\end{equation}
where
\begin{equation}
\chi(W)(y,\bar{y}):=\chi(W,\{id\})(y,\bar{y}).
\end{equation}
\end{rem}
\begin{defn}
Let $W=(a_1,\dots ,a_n;h)$ and $W^*=(a^*_1,\dots ,a^*_n;h^*)$ be reduced
regular weight systems.
Then $W^*$ is called {\it P-dual} to $W$, if
\begin{equation}
\chi(W)(y,\bar{y})=(-1)^n\bar{y}^{\hat{c}_{W^*}}
\chi(W^*,G^*_0)(y,\bar{y}^{-1}),
\end{equation}
and
\begin{equation}
\chi(W^*)(y,\bar{y})=(-1)^n\bar{y}^{\hat{c}_{W}}
\chi(W,G_0)(y,\bar{y}^{-1}),
\end{equation}
where $\hat{c}_W:=1-2\frac{\epsilon_W}{h}$ and $\hat{c}_{W^*}:=1-2
\frac{\epsilon_{W^*}}{h^*}$.
\end{defn}
%
\section{M-Duality for Regular Weight Systems}
In the following section, we assume that $n=3$.
\begin{defn}
The {\it characteristic polynomial} for a regular weight system $W$ is
defined by,
\begin{equation}
\varphi_W(\lambda):=\prod_{i=1}^{\mu_W}(\lambda-{\bf e}[\frac{m_i}{h}]).
\end{equation}
\end{defn}
Let $h/a_i=p_i/q_i$ be the reduced expressions of the rational numbers, i.e.,
\begin{equation}
p_i=\frac{h}{(h,a_i)}~{\rm and }~q_i=\frac{a_i}{(h,a_i)},~~~i=1,2,3.
\end{equation}
Let $p_{ij}:=lcm(p_i,p_j)$, $i,j=1,2,3$ and $p_{123}:=lcm(p_1,p_2,p_3)=h$.
\begin{prop}
$\varphi_W(\lambda)$ has a unique expression:
\begin{equation}
\varphi_W(\lambda)=\prod_{i\in M(W)}(\lambda^i-1)^{e_W(i)},
\end{equation}
where $M(W)$ is a poset $($partial ordered set with respect to the division
relation$)$ given by
\begin{equation}
M(W):=\left\{1,p_1,p_2,p_3,p_{12},p_{23},p_{31},p_{123}=h\right\}.
\end{equation}
We call $e_W(\xi),~\xi\in M(W)$ the {\it cyclotomic exponents}.
\end{prop}
Let $W$ be a regular weight system.
The level of $\xi\in M(W)$ is defined by
\begin{equation}
n(\xi):=\#\{i\in\{1,2,3\}:p_i|\xi\}.
\end{equation}
Two leveled posets $(M(W),n)$ and $(M(W'),n)$ are isomorphic if there exists
an isomorphism $\alpha:M(W)\to M(W')$ of posets with $n(\xi)=n(\alpha(\xi))$.
\begin{thm}[{\rm \cite{sa:2} Theorem 5.2}]\label{th:1}
There are 14 isomorphism classes of the leveled posets $(M(W),n)$ attached to
regular weight systems $W$. The 14 types are exhibited by graphs as follows:
\begin{figure}[hp]
\unitlength=1mm
\begin{picture}(150,45)
\put(0,45){Type}
\put(13,45){$n=0$}
\put(35,45){$n=1$}
\put(65,45){$n=2$}
\put(100,45){$n=3$}
\put(0,10){I.}
\put(15,15){-1}
\put(15.5,11){\circle{6}}
\put(15,9.5){1}

\put(18.5,11){\line(2,1){17.5}}
\put(18.5,11){\line(1,0){15}}
\put(18.5,11){\line(2,-1){16}}

\put(40,25){1}
\put(37,23){\circle{7}}
\put(35.5,23){$p_1$}

\put(40,13){1}
\put(37,11){\circle{7}}
\put(35.5,10){$p_2$}

\put(40,3){1}
\put(37,0){\circle{7}}
\put(35.5,0){$p_3$}

\put(70,26){$-p_1p_2/p_{12}$}
\put(67,23){\oval(8,7)}
\put(65.5,23){$p_{12}$}

\put(70,13){$-p_3p_1/p_{13}$}
\put(67,11){\oval(8,7)}
\put(65.5,10){$p_{13}$}

\put(70,-8){$p_2p_3/p_{23}$}
\put(67,-2){\oval(8,7)}
\put(65.5,-3){$p_{23}$}

\put(41,23){\line(2,-1){22}}
\put(41,23){\line(1,0){22}}

\put(41,11){\line(2,1){22}}
\put(41,11){\line(2,-1){22}}

\put(40,-2){\line(2,1){23}}
\put(40,-2){\line(1,0){23}}

\put(100,15){$p_1p_2p_3/p_{123}$}
\put(100.5,10.5){\oval(9,7)}
\put(97.5,10){$p_{123}$}

\put(71,24){\line(2,-1){25}}
\put(71,11){\line(1,0){24}}
\put(71,-2){\line(2,1){25}}

\end{picture}
\end{figure}
\begin{figure}[hp]
\unitlength=1mm
\begin{picture}(150,30)
\put(0,10){II.}
\put(15,15){-1}
\put(15.5,11){\circle{6}}
\put(15,9.5){1}

\put(18.5,11){\line(2,1){17.5}}
\put(18.5,11){\line(2,-1){16}}

\put(43,23){1}
\put(40,21){\circle{7}}
\put(38,21){$p_1$}

\put(40,3){1}
\put(37,1){\circle{7}}
\put(35.5,1){$p_2$}

\put(70,15){$-p_1p_2/p_{12}$}
\put(67,10){\oval(8,7)}
\put(65.5,9.5){$p_{12}$}

\put(80,-20){$-(p_2-1)/q_{3}$}
\put(72,-12){\oval(20,7)}
\put(65,-12){$p_{3}=p_{23}$}

\put(43,21){\line(2,-1){20}}
\put(40,-0.5){\line(2,1){22}}
\put(40,-0.5){\line(2,-1){21}}

\put(98.5,3){$p_1(p_2-1)p_3/q_3p_{13}$}
\put(108.5,-5){\oval(25,7)}
\put(99.5,-5){$p_{13}=p_{123}$}

\put(71,8){\line(2,-1){24}}
\put(82,-12){\line(2,1){13}}

\end{picture}
\end{figure}
\begin{figure}[hp]
\unitlength=1mm

\begin{picture}(150,30)
\put(0,10){III.}
\put(15,15){-1}
\put(15.5,11){\circle{6}}
\put(15,9.5){1}

\put(18.5,11){\line(2,1){17.5}}
\put(18.5,11){\line(4,-1){40}}

\put(43,23){1}
\put(40,21){\circle{7}}
\put(38,21){$p_1$}

\put(44,21){\line(4,-1){50}}

\put(88,-5){$-(p_2-q_2-q_3)/q_2q_3$}
\put(72,-2){\oval(27,7)}
\put(61,-2){$p_2=p_3=p_{23}$}

\put(83,15){$p_1p_2(p_2-q_2-q_3)/p_{12}q_2p_3$}
\put(108.5,10){\oval(30,7)}
\put(94.5,10){$p_{12}=p_{13}=p_{123}$}

\put(83,2){\line(2,1){11}}

\end{picture}
\end{figure}
\begin{figure}[hp]
\unitlength=1mm
\begin{picture}(150,20)
\put(0,10){IV.}
\put(15,15){-1}
\put(15.5,11){\circle{6}}
\put(15,9.5){1}

\put(18.5,11){\line(1,0){17.5}}

\put(42,15){1}
\put(40,11){\circle{7}}
\put(38,10){$p_1$}

\put(44,11){\line(1,0){10.5}}

\put(56,17){$-(p_1-1)/q_2$}
\put(62,11){\oval(15,7)}
\put(55,10){$p_2=p_{13}$}

\put(84,17){$(p_1-1)(p_2-q_2)/q_2p_3$}
\put(100.5,11){\oval(40,7)}
\put(82.5,10){$p_3=p_{23}=p_{13}=p_{123}$}

\put(70,11){\line(1,0){10}}

\end{picture}
\end{figure}
\begin{figure}[hp]
\unitlength=1mm
\begin{picture}(150,20)
\put(0,10){V.}
\put(15,15){-1}
\put(15.5,11){\circle{6}}
\put(15,9.5){1}

\put(18.5,11){\line(1,0){44}}


\put(90,17){$(\mu_W+1)/h$}
\put(98,11){\oval(70,7)}
\put(65,10){$p_1=p_2=p_3=p_{12}=p_{23}=p_{13}=p_{123}$}
\end{picture}
\end{figure}
\begin{figure}[hp]
\unitlength=1mm
\begin{picture}(150,45)
\put(0,10){VI.}
\put(15,15){-1}
\put(15.5,11){\circle{6}}
\put(15,9.5){1}

\put(18.5,11){\line(2,1){17.5}}
\put(18.5,11){\line(1,0){15}}
\put(18.5,11){\line(2,-1){16}}

\put(40,25){1}
\put(37,23){\circle{7}}
\put(35.5,23){$p_1$}

\put(40,13){1}
\put(37,11){\circle{7}}
\put(35.5,10){$p_2$}

\put(40,3){1}
\put(37,0){\circle{7}}
\put(35.5,0){$p_3$}

\put(70,26){$-p_1p_2/p_{12}$}
\put(67,23){\oval(8,7)}
\put(65.5,23){$p_{12}$}

\put(70,-8){$-p_2p_3/p_{23}$}
\put(67,-2){\oval(8,7)}
\put(65.5,-3){$p_{23}$}


\put(41,23){\line(1,0){22}}
\put(41,11){\line(2,1){22}}
\put(41,11){\line(2,-1){22}}
\put(41,-2){\line(1,0){22}}

\put(95,15){$(p_2-1)p_1p_3/p_{13}$}
\put(106.5,10.5){\oval(20,7)}
\put(97.5,10){$p_{13}=p_{123}$}

\put(71,24){\line(2,-1){25}}

\put(71,-2){\line(2,1){25}}

\end{picture}
\end{figure}
\begin{figure}[hp]
\unitlength=1mm
\begin{picture}(150,40)
\put(0,10){VII.}
\put(15,15){-1}
\put(15.5,11){\circle{6}}
\put(15,9.5){1}

\put(18.5,11){\line(2,1){16}}
\put(18.5,11){\line(1,0){15}}
\put(18.5,11){\line(2,-1){16}}

\put(40,23){1}
\put(37,21){\circle{7}}
\put(35.5,21){$p_1$}

\put(40,13){1}
\put(37,11){\circle{7}}
\put(35.5,10){$p_2$}

\put(40,3){1}
\put(37,0){\circle{7}}
\put(35.5,0){$p_3$}


\put(70,-8){$-p_2p_3/p_{23}$}
\put(67,-2){\oval(8,7)}
\put(65.5,-3){$p_{23}$}


\put(41,21){\line(4,-1){50}}

\put(41,11){\line(2,-1){22}}
\put(41,-2){\line(1,0){22}}

\put(92,15){$p_1(p_2p_3-p_2-p_3)/p_{12}$}
\put(105.9,9){\oval(30,7)}
\put(92,9.2){$p_{12}=p_{13}=p_{123}$}


\put(71,-1.5){\line(2,1){20}}

\end{picture}
\end{figure}
\begin{figure}[hp]
\unitlength=1mm
\begin{picture}(150,40)
\put(0,10){VIII.}
\put(15,15){-1}
\put(15.5,11){\circle{6}}
\put(15,9.5){1}

\put(18.5,11){\line(2,1){16}}
\put(18.5,11){\line(1,0){15}}
\put(18.5,11){\line(2,-1){16}}

\put(40,23){1}
\put(37,21){\circle{7}}
\put(35.5,21){$p_1$}

\put(40,13){1}
\put(37,11){\circle{7}}
\put(35.5,10){$p_2$}

\put(40,3){1}
\put(37,0){\circle{7}}
\put(35.5,0){$p_3$}


\put(41,21){\line(4,-1){40}}

\put(41,11){\line(1,0){40}}
\put(41,0){\line(4,1){40}}

\put(75,18){$(p_1p_2p_3-p_1p_2-p_2p_3-p_3p_1)/p_{12}$}
\put(103,9.5){\oval(43,7)}
\put(83,9.2){$p_{12}=p_{23}=p_{13}=p_{123}$}

\end{picture}
\end{figure}
\begin{figure}[hp]
\unitlength=1mm
\begin{picture}(150,40)
\put(0,10){IX.}
\put(15,15){-1}
\put(15.5,11){\circle{6}}
\put(15,9.5){1}

\put(18.5,11){\line(2,1){16}}
\put(18.5,11){\line(3,-1){16}}

\put(40,24){1}
\put(37,21){\circle{7}}
\put(35.5,21){$p_1$}

\put(40,6){1}
\put(37,4){\circle{7}}
\put(35,4.5){$p_2$}


\put(70,-9){$-(p_2-1)/q_3$}
\put(70,-2){\oval(15,7)}
\put(63,-3){$p_3=p_{23}$}


\put(41,21){\line(4,-1){50}}
\put(41,4){\line(4,-1){21}}

\put(80,17){$p_1(p_2p_3/q_3-p_2-p_3)/q_3)/p_{12}$}
\put(108.1,9){\oval(35,7)}
\put(92,9){$p_{12}=p_{13}=p_{123}$}


\put(76.5,1.4){\line(2,1){14}}

\end{picture}
\end{figure}
\begin{figure}[hp]
\unitlength=1mm
\begin{picture}(150,40)
\put(0,10){X.}
\put(15,15){-1}
\put(15.5,11){\circle{6}}
\put(15,9.5){1}

\put(18.5,11){\line(2,1){15.5}}
\put(18.5,11){\line(3,-1){15.3}}

\put(40,24){1}
\put(37,21){\circle{7}}
\put(35.5,21){$p_1$}

\put(40,6){1}
\put(37,4){\circle{7}}
\put(35,4.5){$p_2$}


\put(65,17){$-p_1p_2/p_{12}$}
\put(66,11){\oval(10,7)}
\put(63,11){$p_{12}$}


\put(41,20.5){\line(3,-1){20}}
\put(41,4.5){\line(3,1){20}}

\put(89,17){$(p_1-1)(p_2-1)/q_3$}
\put(106,11){\oval(40,7)}
\put(88,10){$p_3=p_{23}=p_{13}=p_{123}$}


\put(71.5,11.5){\line(1,0){14}}

\end{picture}
\end{figure}
\begin{figure}[hp]
\unitlength=1mm
\begin{picture}(150,30)
\put(0,10){XI.}
\put(15,15){-1}
\put(15.5,11){\circle{6}}
\put(15,9.5){1}

\put(18.5,11){\line(2,1){15.5}}
\put(18.5,11){\line(3,-1){15.3}}

\put(40,24){1}
\put(37,21){\circle{7}}
\put(35.5,21){$p_1$}

\put(40,6){1}
\put(37,4){\circle{7}}
\put(35,4.5){$p_2$}



\put(41,20.5){\line(4,-1){30}}
\put(41,4.5){\line(4,1){30}}

\put(76,19){$(p_1-1)(p_2-1)/q_3-p_1p_2/p_{12}$}
\put(99,12){\oval(55,7)}
\put(76,12){$p_3=p_{12}=p_{23}=p_{13}=p_{123}$}

\end{picture}
\end{figure}
\begin{figure}[hp]
\unitlength=1mm
\begin{picture}(150,30)
\put(0,10){XII.}
\put(15,15){-1}
\put(15.5,11){\circle{6}}
\put(15,9.5){1}

\put(18.5,11){\line(1,0){15}}

\put(40,13){1}
\put(37,11){\circle{7}}
\put(35.5,10){$p_1$}

\put(66,26){$-(p_1-1)/q_2$}
\put(67,20){\oval(18,7)}
\put(59.5,20){$p_2=p_{12}$}

\put(66,-8){$-(p_1-1)/q_3$}
\put(67,-1){\oval(18,7)}
\put(59.5,-1.5){$p_3=p_{13}$}

\put(41,11){\line(2,1){17}}
\put(41,11){\line(2,-1){17.5}}

\put(95,16){$(p_1-1)p_2p_3/q_2q_3p_{23}$}
\put(106.5,10.5){\oval(20,7)}
\put(97.5,10){$p_{23}=p_{123}$}

\put(76.5,21){\line(2,-1){20}}
\put(76.5,-1){\line(2,1){20}}

\end{picture}
\end{figure}
\begin{figure}[hp]
\unitlength=1mm
\begin{picture}(150,20)
\put(0,10){XIII.}
\put(15,15){-1}
\put(15.5,11){\circle{6}}
\put(15,9.5){1}

\put(18.5,11){\line(1,0){15}}

\put(40,13){1}
\put(37,11){\circle{7}}
\put(35.5,10){$p_1$}


\put(41,11){\line(1,0){21}}

\put(70,16){$(p_1-1)(p_2-q_2-q_3)/q_2q_3$}
\put(92,10.5){\oval(60,7)}
\put(63,10){$p_2=p_3=p_{12}=p_{23}=p_{13}=p_{123}$}


\end{picture}
\end{figure}
\begin{figure}[hp]
\unitlength=1mm
\begin{picture}(150,20)
\put(0,10){XIV.}
\put(15,15){-1}
\put(15.5,11){\circle{6}}
\put(15,9.5){1}

\put(18.5,11){\line(1,0){29}}

\put(40,16){$-(p_{12}-q_1-q_2)/q_1q_2$}
\put(60,10.5){\oval(25,7)}
\put(48,10){$p_1=p_2=p_{12}$}

\put(73,11){\line(1,0){11}}

\put(85,16){$(p_1-q_1)(p_2-q_2)/q_1q_2q_3$}
\put(105,10.5){\oval(42,7)}
\put(85,10){$p_3=p_{23}=p_{13}=p_{123}$}


\end{picture}
\end{figure}
\newpage
where
\begin{enumerate}
\item An element $\xi\in M(W)$ is represented a vertex $\xi$ of a graph.
The vertices are ordered from left to right according to the level $n$.
\item An edge from $\xi$ to $\eta$ is drawn if $\xi|\eta$ and $\xi\ne\eta$
and there does not exist $\zeta$ such that $\zeta\ne\xi$, $\zeta\ne\eta$ and
$\xi|\zeta|\eta$.
\item The number attached near a vertex $\xi$ is the cyclotomic exponent
$e_W(\xi)$.
\end{enumerate}
\end{thm}
\begin{defn}
Let $W$ be a reduced regular weight system. Then
\begin{equation}
mult(W):=e_W(h)
\end{equation}
is called the {\it multiplicity} of $W$.
\end{defn}
Then we can show the following corollary for theorem 2.1.
\begin{cor}{\rm \cite{sa:2}}
Let $W$ be a reduced regular weight system.
Then $mult(W)>0$.
\end{cor}
\begin{rem}
Let $\mu_i$ be the multiplicity of the exponent equal to $i$.
Then we have
\begin{equation}
mult(W)=e_W(h)=\mu_{-1}+\mu_1.
\end{equation}
\end{rem}
\begin{thm}[\cite{sa:2} Theorem 6.2]
Let $W$ be a reduced regular weight system.
Then $mult(W)=1$ if and only if
\begin{equation}
e_W(\xi)=(-1)^{n(\xi)+1},~~~{\rm for~all~}\xi\in M(W).
\end{equation}
\end{thm}
\begin{defn}[\cite{sa:2} Definition 7.5]
Let $W$ and $W^*$ be reduced regular weight systems.
$W^*$ is said to be {\it M-dual} to $W$, if
\begin{equation}
\prod_{i\in M(W^*)}(\lambda^i-1)^{e_{W^*}(i)}=\prod_{i\in M(W)}(\lambda^i-1)^
{-e_W(h/i)}.
\end{equation}
\end{defn}
\begin{rem}(\cite{sa:2} Assertion 7.7)
Let $W$ be a reduced regular weight system.
If $W$ admits a M-dual $W^*$, then
\begin{enumerate}
\item $\mu_0=0$,
\item $mult(W)=1$,
\item $(M(W),n)$ is either of types I, II, III, IV or V.
\end{enumerate}
\end{rem}
Let $W$ be a reduced regular weight system of type V and $mult(W)=1$.
Then $W$ is given by
\begin{equation}
W=(lm-m+1,mk-k+1,kl-l+1;h),
\end{equation}
for $k,l,m\in\Z_{\ge 0}$ such that $h_W=klm+1$ and $(lm-m+1,h_W)=1$.
The characteristic polynomial for $W$ is given by
\begin{equation}
\varphi_{W}=\frac{(\lambda^{h_W}-1)}{(\lambda-1)}.
\end{equation}
This depends only on $h_W$.
So we define the M-duality of type V as follows:
\begin{defn}
Let $W$ and $W^*$ be reduced regular weight systems of type V with $mult(W)=1$
and $mult(W^*)=1$.
$W^*$ is said to be M-dual to $W=(lm-m+1,mk-k+1,kl-l+1;klm+1)$,
if $W^*=(lm-l+1,mk-m+1,kl-k+1;klm+1)$.
\end{defn}
\begin{note}
Let $W$ be M-dual to $W^*$.
We give a list of $W$ and $W^*$ for each type of $M(W)$.
\bigskip
\noindent
{\bf I.}
$M(W)$ and $M(W^*)$ are of type I with $mult(W)=1$ and $mult(W^*)=1$,
if and only if $(p_i,p_j)=1$, $i=1,2,3$.
\begin{equation}
W=W^*=(p_2p_3,p_3p_1,p_1p_2;p_1p_2p_3).
\end{equation}
\bigskip
\noindent
{\bf II.} $M(W)$ and $M(W^*)$ are of type II with $mult(W)=1$ and
$mult(W^*)=1$, if and only if $p_2\ne p_3$, $p_2|p_3$, $(p_1,p_3)=1$,
$(p_2-1,p_3)=1$ and $(p_3/p_2-1,p_3)=1$.
\begin{equation}
W=(p_3,\frac{p_1p_3}{p_2},(p_2-1)p_1;p_1p_3),
\end{equation}
\begin{equation}
W^*=(p_3,p_1p_2,(\frac{p_3}{p_2}-1)p_1;p_1p_3).
\end{equation}
\bigskip
\noindent
{\bf III.} $M(W)$ and $M(W^*)$ are of type III with $mult(W)=1$ and
$mult(W^*)=1$, if and only if $p_2=p_3$, $(p_1,p_2)=1$.
\begin{equation}
W=W^*=(p_2,p_1q_2,p_1q_3;p_1p_2),
\end{equation}
for any integer $q_2$ and $q_3$ such that $p_2+1=(q_2+1)(q_3+1)$ and
$(q_2,q_3)=1$.
\bigskip
\noindent
{\bf IV.} $M(W)$ and $M(W^*)$ are of type IV with $mult(W)=1$ and
$mult(W^*)=1$, if and only if $p_1\ne p_2\ne p_3$, $p_1|p_3$, $p_2|p_3$,
$(p_1-1,p_2)=1$, $(p_2-p_1+1,p_3)=1$, $(p_3/p_2-1,p_3/p_1)=1$,
$(p_3/p_1-p_3/p_2+1,p_3)=1$.
\begin{equation}
W=(\frac{p_3}{p_1},(p_1-1)\frac{p_3}{p_2},p_2-p_1+1;p_3),
\end{equation}
\begin{equation}
W^*=(p_2,(\frac{p_3}{p_2}-1)p_1,\frac{p_3}{p_1}-\frac{p_3}{p_2}+1;p_3).
\end{equation}
\bigskip
\noindent
{\bf V.} $M(W)$ and $M(W^*)$ are of type V with $mult(W)=1$ and $mult(W^*)=1$,
if and only if $p_1=p_2=p_3=h$.
\begin{equation}
W=(lm-m+1,mk-k+1,kl-l+1;h),
\end{equation}
\begin{equation}
W^*=(lm-l+1,mk-m+1,kl-k+1;h),
\end{equation}
for any positive integers $k,l,m$ such that $h=klm+1$, $(lm-m+1,h)=1$ and
$(lm-l+1,h)=1$.
\end{note}
%
\section{Relation between M-Duality and P-Duality}
In this section, we will prove the following theorem:
\begin{thm}
Let $W$ and $W^*$ be the reduced regular weight system.
Then $W^*$ is P-dual to $W$, if and only if $W^*$ is M-dual to $W$.
\end{thm}
\begin{lem}
Let $W$ and $W^*$ be a reduced regular weight systems.
If $W^*$ is P-dual to $W$, then the multiplicity of the exponent equal to $0$
is $0$, i.e., $\mu_0=0$ $(\mu^*_0=0)$.
\end{lem}
\begin{pf}
Since $W^*$ is P-dual to $W$, we have
\begin{equation}
\chi(W^*)(y,\bar{y})=(-1)^3\bar{y}^{\hat{c}_W}
\chi(W,G_0)(y,\bar{y}^{-1}).
\label{eq:1}
\end{equation}
The right-hand side of \eqref{eq:1} is
\begin{equation}
\begin{split}
(-1)^3\bar{y}^{\hat{c}_W}\chi(W,G_0)(y,\bar{y}^{-1}) &= (-1)^3
\bar{y}^{\hat{c}_W}\chi_0(W,G_0)(y,\bar{y}^{-1})\\
& +(-1)^3\bar{y}^{\hat{c}_W}\sum_{l_2\in L_2}\chi_{l_2}(W,G_0)
(y,\bar{y}^{-1})\\ & +(-1)^3\bar{y}^{\hat{c}_W}\sum_{l_1\in L_1}\chi_{l_1}
(W,G_0)(y,\bar{y}^{-1})\\ &+(-1)^n\bar{y}^{\hat{c}_W}\sum_{l_0\in L_0}
\chi_{l_0}(W,G_0)(y,\bar{y}^{-1}),
\end{split}
\end{equation}
where
\begin{equation}
L_2:=\{0<l<h_{W}~|~\xi|l~{\rm for}~\xi\in M(W)~s.t.~n(\xi)=2\},
\end{equation}
\begin{equation}
L_1:=\{0<l<h_{W}~|~\xi|l~{\rm for}~\xi\in M(W)~s.t.~n(\xi)=1\}
\backslash L_2,
\end{equation}
and
\begin{equation}
L_0:=\{0<l<h_{W}~|~ \xi\nmid l~{\rm for~all}~\xi\in M(W)\}.
\end{equation}
We have
\begin{equation}
(-1)^3\bar{y}^{\hat{c}_W}\chi_0(W,G_0)(y,\bar{y}^{-1})=-\mu_0
y^{-\frac{\epsilon_W}{h_W}}\bar{y}^{1-\frac{\epsilon_W}{h_W}}
-\mu_{h_W}y^{1-\frac{\epsilon_W}{h_W}}\bar{y}^{-\frac{\epsilon_W}{h_W}},
\end{equation}
\begin{equation}
(-1)^3\bar{y}^{\hat{c}_W}\sum_{l_2\in L_2}\chi_{l_2}(W,G_0)
(y,\bar{y}^{-1})=\sum_{l_2\in L_2}\nu_{l_2}(y\bar{y})^{f_{l_2}},
\end{equation}
where $\nu_{l_2}$ is a non--negative integer and
\begin{equation}
f_{l_2}:=1-\left(\frac{q_i}{p_i}l_2-\left[\frac{q_i}{p_i}l_2
\right]\right)-\frac{\epsilon_W}{h_W},~~~p_i\nmid l_2,~l_2\in L_2,
\end{equation}
\begin{equation}
(-1)^3\bar{y}^{\hat{c}_W}\sum_{l_1\in L_1}\chi_{l_1}(W,G_0)
(y,\bar{y}^{-1})=0,
\end{equation}
and
\begin{equation}
(-1)^3\bar{y}^{\hat{c}_W}\sum_{l_0\in L_0}\chi_{l_0}(W,G_0)
(y,\bar{y}^{-1})=\sum_{l_0\in L_0}(y\bar{y})^{f_{l_0}},
\end{equation}
where
\begin{equation}
\begin{split}
f_{l_0}:=& 2-\sum_{i=1}^3\left(\frac{q_i}{p_i}l_0-\left[\frac{q_i}{p_i}
l_0\right]\right)-\frac{\epsilon_W}{h_W}\\
=& 2-l_0-\frac{\epsilon_W}{h_W}l_0+\sum_{i=1}^3\left[\frac{q_i}{p_i}
l_0\right]-\frac{\epsilon_W}{h_W},~~~l_0\in L_0.
\end{split}
\end{equation}
Since the left-hand side of \eqref{eq:1} can not contain monomials as
$y^a\bar{y}^b$, $a\ne b$.
We have $\mu_0=\mu_{h_W}=0$.
\end{pf}\qed
\begin{lem}
Let $W$ and $W^*$ be a reduced regular weight system.
If $W^*$ is P-dual to $W$, then
\begin{equation}
\frac{\epsilon_{W}}{h_{W}}=\frac{\epsilon_{W^*}}{h_{W^*}}.
\end{equation}
\end{lem}
\begin{pf}
The maximal exponent of $(y\bar{y})$ in the left-hand side of \eqref{eq:1} is
$1-2\frac{\epsilon_{W^*}}{h_{W^*}}$ and the minimal one is $0$.
So we will show the maximal exponent of $(y\bar{y})$ in the right-hand side of
\eqref{eq:1} is $1-2\frac{\epsilon_W}{h_W}$.
It is obvious that if $l_i\in L_i$, then $h-l_i\in L_i$ (for $i=0,2$).
Thus we have the duality property of exponents
\begin{equation}
f_{l_i}+f_{h_W-l_i}=1-2\frac{\epsilon_W}{h_W},~~~i=0,2.
\end{equation}
Note that $1$ is always contained in $L_0$.
We have $f_1=1-2\frac{\epsilon_W}{h_W}$ and $f_{h_W-1}=0$.
Since $0$ is minimal, using the duality property above, $f_1$ is maximal.
\end{pf}\qed
\begin{lem}
Let $W$ and $W^*$ be a reduced regular weight system.
If $W^*$ is P-dual to $W$, then
$(\epsilon_W,h_W)=1$ $((\epsilon_{W^*},h_{W^*})=1)$.
\end{lem}
\begin{pf}
Suppose that there exists a prime number $k$ such that $k|(\epsilon_W,h_W)$.
This implies that $k|\xi$ for all $\xi\in M(W)$, $n(\xi)=2$.
(Since $k|(h_W=p_{123})$, we can assume that $k|p_1$, $k\nmid p_2$ and
$k\nmid p_3$. This is impossible since $k|\epsilon_W$ implies that
\begin{equation}
k\left|\left(q_1\frac{p_{23}}{(p_1,p_{23})}+q_2\frac{p_{31}}{(p_2,p_{31})}+
q_3\frac{p_{12}}{(p_3,p_{12})}\right)\right.,
\end{equation}
this means $k|p_{23}$.)
If $k|\xi$ for all $\xi\in M(W)$, $n(\xi)=2$, it is clear that
\begin{equation}
f_{l_i}\cdot\frac{h_W}{k}\in\Z,~~~l_i\in L_i,~~i=0,2.
\end{equation}
This contradicts that $mult(W^*)>0$ (i.e. the existence of exponents prime to
$h_{W^*}$).
\end{pf}\qed
\begin{cor}
Let $W$ and $W^*$ be a reduced regular weight system.
If $W^*$ is P-dual to $W$, then $h_W=h_{W^*}$.
\end{cor}
When $W^*$ is P-dual to $W$, we will denote $h_W=h_{W^*}$ by $h$ and
$\epsilon_W=\epsilon_{W^*}$ by $\epsilon$.
\begin{lem}
Let $W$ and $W^*$ be a reduced regular weight system.
If $W^*$ is P-dual to $W$, then $mult(W)=1$ $(mult(W^*)=1)$.
\end{lem}
\begin{pf}
Since $(\epsilon,h)=1$, it is clear that there does not exist
$l_i\in L_i$, $i=0,2$ such that $f_{l_i}+\epsilon/h=-\epsilon/h$.
\end{pf}\qed
Now we give a proof of the theorem 3.1.
\begin{pf}
We will prove by the classification of the poset $M(W)$ (theorem \ref{th:1})
and the duality in \cite{bh:1}.
Note that for regular reduced weight systems $W$ and $W'$,
\begin{equation}
\chi(W)(y,\bar{y})=\chi(W')(y,\bar{y}),
\end{equation}
if and only if $W=W'$.
Let $W$ be a reduced regular weight system with $mult(W)=1$.
\bigskip
\noindent
{\bf I.}
$W$ is given by
\begin{equation}
W=(p_2p_3,p_3p_1,p_1p_2;p_1p_2p_3).
\end{equation}
We can take the weighted homogeneous polynomial $F_W(x,y,z)$ associated to $W$:
\begin{equation}
F_W(x,y,z):=x^{p_1}+y^{p_2}+z^{p_3}.
\end{equation}
We see that the group of phase symmetries $G_{F_W}$ is isomorphic to $G_0$.
The transposed polynomial \cite{bh:1} $F^*_W(x_*,y_*,z_*)$ is given by
\begin{equation}
F^*_W(x_*,y_*,z_*)=x_*^{p_1}+y_*^{p_2}+z_*^{p_3},
\end{equation}
and the group of phase symmetries $G_{F^*_W}$ is obviously isomorphic to
$\Z_{h}$.
Then it is well--known \cite{ge:1}\cite{gq:1} that
\begin{equation}
\chi(W,\{id\})(y,\bar{y})=(-1)^3\bar{y}^{\frac{h-2\epsilon}{h}}\chi(W,G_{F_W})
(y,\bar{y}^{-1}).
\end{equation}
So $W^*=W$ and $W^*$ is P-dual to $W$ if and only if $W^*$ is M-dual to $W$.
\bigskip
\noindent
{\bf II.} $W$ is given by
\begin{equation}
W=(p_3,\frac{p_1p_3}{p_2},(p_2-1)p_1;p_1p_3).
\end{equation}
We assume that $(\epsilon,h)=1$.
We can take the weighted homogeneous polynomial $F_W(x,y,z)$ associated to $W$:
\begin{equation}
F_W(x,y,z):=x^{p_1}+y^{p_2}+yz^{\frac{p_3}{p_2}}.
\end{equation}
We see that the group of phase symmetries $G_{F_W}$ is isomorphic to $G_0$.
The transposed polynomial $F^*_W(x_*,y_*,z_*)$ is given by
\begin{equation}
F^*_W(x_*,y_*,z_*)=x_*^{p_1}+y_*^{\frac{p_3}{p_2}}+y_*z_*^{p_2},
\end{equation}
and the group of phase symmetries $G_{F^*_W}$ is isomorphic to $\Z_{h}$.
Then since $(p_1,p_3)=1$, $(p_2-1,p_3)=1$ and $(\epsilon,h)=1$,
using the duality of {\it the linear chain} in \cite{bh:1}, we see that
\begin{equation}
\chi(W^*,\{id\})(y,\bar{y})=(-1)^3\bar{y}^{\frac{h-2\epsilon}{h}}
\chi(W,G_{F_W})(y,\bar{y}^{-1}),
\end{equation}
and
\begin{equation}
\chi(W,\{id\})(y,\bar{y})=(-1)^3\bar{y}^{\frac{h-2\epsilon}{h}}
\chi(W^*,G_{F^*_W})(y,\bar{y}^{-1}),
\end{equation}
where $W^*$ is the regular weight system associated to $F^*_W$ and given by
\begin{equation}
W^*=(p_3,p_1p_2,(\frac{p_3}{p_2}-1)p_1;p_1p_3).
\end{equation}
So $W^*$ is P-dual to $W$ if and only if $W^*$ is M-dual to $W$.
\bigskip
\noindent
{\bf III.} $W$ is given by
\begin{equation}
W=(p_2,p_1q_2,p_1q_3;p_1p_2),
\end{equation}
for any integer $q_2$ and $q_3$ such that $p_2+1=(q_2+1)(q_3+1)$ and
$(q_2,q_3)=1$.
We can take the weighted homogeneous polynomial $F_W(x,y,z)$ associated to $W$:
\begin{equation}
F_W(x,y,z):=x^{p_1}+y^{q_3+1}z+yz^{q_2+1}.
\end{equation}
The group of phase symmetries $G_{F_W}$ is isomorphic to $\Z_{h}$.
The transposed polynomial $F^*_W(x_*,y_*,z_*)$ is given by
\begin{equation}
F^*_W(x_*,y_*,z_*)=x_*^{p_1}+y_*^{q_3+1}z_*+y_*z_*^{q_2+1},
\end{equation}
and the group of phase symmetries $G_{F^*_W}$ is isomorphic to
$\Z_{h}$.
Then since $(p_1,p_2)=1$, $(q_2,p_2)=1$ and $(q_3,p_2)=1$, using the duality
of {\it the loop} in \cite{bh:1}, we see that
\begin{equation}
\chi(W,\{id\})(y,\bar{y})=(-1)^3\bar{y}^{\frac{h-2\epsilon}{h}}
\chi(W,G_{F_W})(y,\bar{y}^{-1}),
\end{equation}
So $W^*=W$ and $W^*$ is P-dual to $W$ if and only if $W^*$ is M-dual to $W$.
\bigskip
\noindent
{\bf IV.}
$W$ is given by
\begin{equation}
W=(\frac{p_3}{p_1},(p_1-1)\frac{p_3}{p_2},p_2-p_1+1;p_3).
\end{equation}
We assume that $(\epsilon,h)=1$.
We can take the weighted homogeneous polynomial $F_W(x,y,z)$ associated to $W$:
\begin{equation}
F_W(x,y,z):=x^{p_1}+xy^{\frac{p_2}{p_1}}+yz^{\frac{p_3}{p_2}}.
\end{equation}
The group of phase symmetries $G_{F_W}$ is isomorphic to $\Z_{h}$.
The transposed polynomial $F^*_W(x_*,y_*,z_*)$ is given by
\begin{equation}
F^*_W(x_*,y_*,z_*)=x_*^{\frac{p_3}{p_2}}+x_*y_*^{\frac{p_2}{p_1}}+y_*z_*^{p_1},
\end{equation}
and the group of phase symmetries $G_{F^*_W}$ is isomorphic to
$\Z_{h}$.
Then since $(p_2-p_1+1,p_3)=1$ and $(\epsilon,h)=1$, using the duality of
{\it the linear chain} in \cite{bh:1}, we see that
\begin{equation}
\chi(W^*,\{id\})(y,\bar{y})=(-1)^3\bar{y}^{\frac{h-2\epsilon}{h}}
\chi(W,G_{F_W})(y,\bar{y}^{-1}),
\end{equation}
and
\begin{equation}
\chi(W,\{id\})(y,\bar{y})=(-1)^3\bar{y}^{\frac{h-2\epsilon}{h}}
\chi(W^*,G_{F^*_W})(y,\bar{y}^{-1}),
\end{equation}
where $W^*$ is the regular weight system associated to $F^*_W$ and given by
\begin{equation}
W^*=(p_2,(\frac{p_3}{p_2}-1)p_1,\frac{p_3}{p_1}-\frac{p_3}{p_2}+1;p_3).
\end{equation}
So $W^*$ is P-dual to $W$ if and only if $W^*$ is M-dual to $W$.
\bigskip
\noindent
{\bf V.} $W$ is given by
\begin{equation}
W=(lm-m+1,kl-l+1,mk-k+1;h),
\end{equation}
for any positive integers $k,l,m$ such that $h=klm+1$ and $(lm-m+1,h)=1$.
We assume that $(\epsilon,h)=1$.
We can take the weighted homogeneous polynomial $F_W(x,y,z)$ associated to $W$:
\begin{equation}
F_W(x,y,z):=zx^{k}+xy^{m}+yz^{l}.
\end{equation}
The group of phase symmetries $G_{F_W}$ is isomorphic to $\Z_{h}$.
The transposed polynomial $F^*_W(x_*,y_*,z_*)$ is given by
\begin{equation}
F^*_W(x_*,y_*,z_*)=z_*x_*^{k}+x_*y_*^{l}+y_*z_*^{m},
\end{equation}
and the group of phase symmetries $G_{F^*_W}$ is isomorphic to
$\Z_{h}$.
Then since $(lm-m+1,h)=1$ and $(\epsilon,h)=1$, using the duality of
{\it the loop} in \cite{bh:1}, we see that
\begin{equation}
\chi(W^*,\{id\})(y,\bar{y})=(-1)^3\bar{y}^{\frac{h-2\epsilon}{h}}
\chi(W,G_{F_W})(y,\bar{y}^{-1}),
\end{equation}
and
\begin{equation}
\chi(W,\{id\})(y,\bar{y})=(-1)^3\bar{y}^{\frac{h-2\epsilon}{h}}
\chi(W^*,G_{F^*_W})(y,\bar{y}^{-1}),
\end{equation}
where $W^*$ is the regular weight system associated to $F^*_W$ and given by
\begin{equation}
W^*=(lm-l+1,mk-m+1,kl-k+1;h).
\end{equation}
So $W^*$ is P-dual to $W$ if and only if $W^*$ is M-dual to $W$.
\end{pf}\qed
\begin{rem}
The characteristic polynomials $\varphi_{W^*}(\lambda)$ are calculated
as follows.
We define the characteristic polynomials $\varphi_{L_0}(\lambda)$ and
$\varphi_{L_2}(\lambda)$:
\begin{equation}
\varphi_{L_0}(\lambda):=\prod_{l_0\in L_0}
\left(\lambda-{\bf e}(f_{l_0}+\frac{\epsilon}{h})\right),
\end{equation}
\begin{equation}
\varphi_{L_2}(\lambda):=\prod_{l_2\in L_2}
\left(\lambda-{\bf e}(f_{l_2}+\frac{\epsilon}{h})\right)^{\nu_{l_2}}.
\end{equation}
\bigskip
\noindent
{\bf I.}
$\varphi_{L_0}(\lambda)$ is given by,
\begin{equation}
\varphi_{L_0}(\lambda)=\frac{(\lambda^{p_1p_2p_3}-1)(\lambda^{p_1}-1)
(\lambda^{p_2}-1)(\lambda^{p_3}-1)}{(\lambda^{p_1p_2}-1)(\lambda^{p_2p_3}-1)
(\lambda^{p_3p_1}-1)
(\lambda -1)},
\end{equation}
since $(\epsilon,h)=1$.
The multiplicity of the exponents $\nu_{l_2}=0$ for $l_2\in L_2$, since
\begin{multline}
\nu_{l_2}=\frac{1}{p_1p_2p_3}( (1-p_1)(1-p_2)p_3+(1-p_1)(p_2p_3-p_3)\\
+(1-p_2)(p_1p_3-p_3)+(p_1p_2p_3-p_2p_3-p_1p_3+p_3) )=0.
\end{multline}
So we have
\begin{equation}
\begin{split}
\varphi_{W^*}(\lambda)&=\varphi_{L_0}(\lambda)\\
&=\frac{(\lambda^h-1)(\lambda^{p_1}-1)(\lambda^{p_2}-1)
(\lambda^{p_3}-1)}{(\lambda^{p_1p_2}-1)(\lambda^{p_2p_3}-1)(\lambda^{p_3p_1}-1)
(\lambda -1)}.
\end{split}
\end{equation}
\bigskip
\noindent
{\bf II.}
\begin{equation}
\varphi_{L_0}(\lambda)=\frac{(\lambda^{p_1p_3}-1)(\lambda^{\frac{p_3}{p_2}}-1)}
{(\lambda^{p_3}-1)(\lambda^{\frac{p_1p_3}{p_2}}-1)}.
\end{equation}
The multiplicity of the exponents $\nu_{l_2}=0$ for $l_2\in L_2$, $p_1p_2l_2$
and $\nu_{l_2}=1$ for $l_2\in L_2$, $p_3|l_2$, since
\begin{multline}
\nu_{l_2}=\frac{1}{p_1p_3}((1-p_1)(1-p_2)\frac{p_3}{p_2}+(1-p_1)
(p_3-\frac{p_3}{p_2})+(1-p_2)(\frac{p_1p_3}{p_2}-\frac{p_3}{p_2})\\
+(p_1p_3+\frac{p_3}{p_2}-p_3-\frac{p_1p_3}{p_2}))=0,
~~~l_2\in L_2,~p_1p_2|l_2,
\end{multline}
\begin{equation}
\begin{split}
\nu_{l_2}&=\frac{1}{p_1p_3}\left((1-p_2)(1-\frac{p_3}{p_2-1})p_1+(1-p_2)
(\frac{p_1p_3}{p_2}-p_1)+(p_1p_3-\frac{p_1p_3}{p_2})\right)\\
&=1,~~~l_2\in L_2,~p_3|l_2.
\end{split}
\end{equation}
$\varphi_{L_2}(\lambda)$ is given by:
\begin{equation}
\varphi_{L_2}(\lambda)=\frac{(\lambda^{p_1}-1)}{(\lambda-1)}.
\end{equation}
Then
\begin{equation}
\begin{split}
\varphi_{W^*}(\lambda)&=\varphi_{L_0}(\lambda)\cdot\varphi_{L_2}(\lambda)\\
&=\frac{(\lambda^{p_1p_3}-1)(\lambda^{p_1}-1)(\lambda^{
\frac{p_3}{p_2}}-1)}{(\lambda^{p_3}-1)(\lambda^{\frac{p_1p_3}{p_2}}-1)
(\lambda -1)}.
\end{split}
\end{equation}
\bigskip
\noindent
{\bf III.}
\begin{equation}
\varphi_{L_0}(\lambda)=\frac{(\lambda^{p_1p_2}-1)(\lambda-1)}
{(\lambda^{p_2}-1)(\lambda^{p_1}-1)},
\end{equation}
\begin{equation}
\begin{split}
\nu_{l_2}&=\frac{1}{p_1p_2}\left((1-\frac{p_1}{q_1})(1-\frac{p_2}{q_2})p_1+
(p_1p_2-p_1)\right)\\
&=2,~~~l_2\in L_2,
\end{split}
\end{equation}
and
\begin{equation}
\varphi_{L_2}(\lambda)=\frac{(\lambda^{p_1}-1)^2}{(\lambda-1)^2}.
\end{equation}
Thus
\begin{equation}
\begin{split}
\varphi_{W^*}(\lambda)&=\varphi_{L_0}(\lambda)\cdot\varphi_{L_2}(\lambda)\\
&=\frac{(\lambda^{p_1p_2}-1)(\lambda^{p_1}-1)}
{(\lambda^{p_2}-1)(\lambda -1)},
\end{split}
\end{equation}
\bigskip
\noindent
{\bf IV.}
\begin{equation}
\varphi_{L_0}(\lambda)=\frac{(\lambda^{p_3}-1)}
{(\lambda^{\frac{p_3}{p_1}}-1)},
\end{equation}
\begin{equation}
\begin{split}
\nu_{l_2}&=\frac{1}{p_3}\left((1-\frac{p_2}{(p_1-1)})(1-p_1)\frac{p_3}{p_1}
+(1-p_1)(\frac{p_3}{p_1}-\frac{p_3}{p_2})+(p_3-\frac{p_3}{p_1})\right)\\
&=2,~~~l_2\in L_2,
\end{split}
\end{equation}
and
\begin{equation}
\varphi_{L_2}(\lambda)=\frac{(\lambda^{\frac{p_3}{p_2}}-1)}{(\lambda-1)}.
\end{equation}
Then
\begin{equation}
\begin{split}
\varphi_{W^*}(\lambda)&=\varphi_{L_0}(\lambda)\cdot\varphi_{L_2}(\lambda)\\
&=\frac{(\lambda^{p_3}-1)(\lambda^{\frac{p_3}{p_2}}-1)}
{(\lambda^{\frac{p_3}{p_1}}-1)(\lambda -1)}.
\end{split}
\end{equation}
\bigskip
\noindent
{\bf V.}
\begin{equation}
\varphi_{L_0}(\lambda)=\frac{(\lambda^{klm+1}-1)}
{(\lambda-1)}.
\end{equation}
Since $L_2=\emptyset$, we have
\begin{equation}
\begin{split}
\varphi_{W^*}(\lambda)&=\varphi_{L_0}(\lambda)\\
&=\frac{(\lambda^{klm+1}-1)}
{(\lambda -1)}.
\end{split}
\end{equation}
\end{rem}
%

%

\newpage
\begin{thebibliography}{[BH]}
%
\bibitem[A]{ar:1}
	V.~I.~Arnold,
	Proc.~Internat.~Congress Math., Vancouver, I (1974) 19.
\bibitem[BH1]{bh:1}
	P.~Berglund, M.~Henningson,
	Nucl.~Phys. {\bf B433} (1994) 311.
\bibitem[BH2]{bh:2}
	P.~Berglund, M.~Henningson, in \cite{gy:1}.
\bibitem[DN]{dn:1}
	I.~V.~Dolgachev, V.~V.~Nikulin,
	Seventh all--union topological conference, Minsk 1977.
\bibitem[G]{ge:1}
	D.~Gepner,
	Nucl.~Phys. {\bf B296} (1988) 757.
\bibitem[GQ]{gq:1}
	D.~Gepner, Z.~Qiu,
	Nucl.~Phys. {\bf B285} (1987) 423.
\bibitem[GY]{gy:1}
	B.~Greene, S.~-T.~Yau, ed.,
	{\it Mirror Symmetry II},
	AMS/IP, 1997.
\bibitem[IV]{iv:1}
	K.~Intriligator, C.~Vafa,
	Nucl.~Phys. {\bf B339} (1990) 95.
\bibitem[KY]{ka:1}
	T.~Kawai, S.-K.~Yang,
	Prog.~Theor.~Phys.~Suppl. No.{\bf 118} (1995) 277.
\bibitem[P]{pi:1}
	H.~Pinkham,
	C.~R.~Acad.~Sc.~Paris {\bf 284 A} (1977) 615.
\bibitem[S1]{sa:1}
	K.~Saito,
	Publ.~RIMS,~Kyoto University {\bf 19} (1983) 1231.
\bibitem[S2]{sa:2}
	K.~Saito,
	{\it Duality for Regular Systems of Weights},
	preprint.
\bibitem[T]{ta:1}
	A.~Takahashi,
	{\it Primitinve Forms and Topological Landau--Ginzburg Model Coupled to
	Gravity}, to appear in master thesis.
%
\end{thebibliography}
\end{document}